\documentclass[a4paper,11pt]{article}
\usepackage{amsmath,amsfonts,amssymb,amscd,amsthm,epsfig}
\usepackage[mathscr]{euscript}
\usepackage{hyperref}
\usepackage{color}
\newtheorem{thm}{Theorem}[section]
\newtheorem{lemma}[thm]{Lemma}
\newtheorem{prop}[thm]{Proposition}
\newtheorem{corr}[thm]{Corollary}
\theoremstyle{definition}
\newtheorem*{nota}{Notation}
\newtheorem*{exmple*}{Example}
\newtheorem{rmks}[thm]{Remarks}
\newtheorem{dfn}[thm]{Definition}

\newtheorem{propr}[thm]{Property}
\newtheorem{proprs}[thm]{Properties}

\newtheorem{constr}[thm]{Construction}
\newtheorem*{thm*}{Theorem}
\theoremstyle{remark}
 \newtheorem*{claim}{\textbf{Claim}}
\newtheorem*{rmq}{\textit{Remark}}
\newtheorem{rmk}[thm]{\textit{Remark}}
\renewcommand{\proof}{\noindent\textit{Proof}\/: \,\,}
\newcounter{lijstc}
{\end{list}}
{\end{list}}
\newenvironment{lijstb}{\vspace{-0.3ex}\begin{list}{\rm \Alph{lijstc})}{\usecounter{lijstc}%
\setlength{\leftmargin}{1ex}
\setlength{\rightmargin}{\leftmargin}}}%
{\end{list}}
\newcommand{\C}{{\mathbb{C}}}
\newcommand{\D}{{\mathbb{D}}}
\newcommand{\Q}{{\mathbb{Q}}}
\newcommand{\R}{{\mathbb{R}}}

\newcommand{\Z}{{\mathbb{Z}}}
\newcommand{\bV}{{\mathbb{V}}}

\newcommand{\bH}{{\mathbb{H}}}
\newcommand{\bP}{{\mathbb{P}}}
\newcommand\CC{{\mathcal C}}

\newcommand\MM{{\mathcal M}}
\newcommand{\II}{{\mathcal I}}
\newcommand{\comp}{\raise1pt\hbox{{$\scriptscriptstyle\circ$}}}
 \def\lset{\{} % for {
\def\rset{\}} % for }
\def\set#1{\lset#1\rset}
\renewcommand\setminus{-}
\newcommand\Tr{{}^{\mathsf{T}}\kern-0.9pt} %transpose
\def\id{\mathop{\rm id}\nolimits}
\newcommand\Hom{\mathop{\rm Hom}\nolimits}
\newcommand\mhc{\mathsf{MHC}}
\newcommand\mhm{\mathsf{MHM}}
\newcommand\mhw{\mathsf{MHW}}
\newcommand\perv{\mathop{\rm Perv}\nolimits}
\newcommand\rat{\mbox{\rm rat}}
\newcommand\Gr{\mathop{\rm Gr}\nolimits}
\newcommand\ic{\mathop{\II\CC}\nolimits}
\newcommand\supp{\mathop{\rm supp}}
\newcommand\DR{\mathop{\rm DR}}
\newcommand\Dec{\mathop{\rm Dec}}
\newcommand\Ker{\mathop{\rm Ker}}
\newcommand\cone{\mathop{\rm Cone}}
\newenvironment{diagram}{
\begin{matrix}}{\end{matrix}}
 \def\arrow(#1,#2)\dir(#3,#4)\long#5{\put(#1,#2){\vector(#3,#4){#5}}}
\def\mapright#1{\mathop{\vbox{\ialign{
                                ##\crcr
    ${\scriptstyle\hfil\;\;#1\;\;\hfil}$\crcr
 \noalign{\kern2pt\nointerlineskip}
    \rightarrowfill\crcr}}\;}}

\def\mapleft#1{\mathop{\vbox{\ialign{
                                ##\crcr
    ${\scriptstyle\hfil\;\;#1\;\;\hfil}$\crcr
 \noalign{\kern2pt\nointerlineskip}
    \leftarrowfill\crcr}}\;}}
\newcommand\Vline[3]{\llap{$\scriptstyle #1$}
\left\Vert\vbox to#3{}\right.\rlap{$\scriptstyle #2$}}
\def\into{\hookrightarrow}

\date{December 10,  2008}
\title{Lowest Weights in Cohomology of Variations of Hodge Structure}
\author{Chris PETERS \\
Department of Mathematics,  University of Grenoble\\
UMR 5582 CNRS-UJF, 38402-Saint-Martin d'H{\`e}res, France\\ 
{\tt chris.peters@ujf-grenoble.fr} \\
and\\
Morihiko SAITO \\
RIMS Kyoto University, Kyoto 606-8502 Japan\\ 
{\tt msaito@kurims.kyoto-u.ac.jp}}
\begin{document}
\maketitle

\begin{abstract} Let $X$ be an irreducible complex analytic space with $j:U\into X$ an immersion of a smooth Zariski open subset, and let $\bV$ be a variation of Hodge structure of weight $n$ over $U$. Assume $X$ is compact K\"ahler. Then provided the local monodromy operators at infinity are quasi-unipotent,  $IH^k(X, \bV)$ is known to carry a pure Hodge structure of weight $k+n$, while $H^k(U,\bV)$ carries a mixed Hodge structure of weight $\ge k+n$. In this note it is shown that the image of the natural map $IH^k(X,\bV) \to H^k(U,\bV)$ is the lowest weight part of this mixed Hodge structure. In the algebraic case this easily follows from the formalism of mixed sheaves, but the analytic case is rather complicated, in particular  when the complement $X-U$ is not a hypersurface.

\end{abstract}

\section*{Introduction}
For a compact K\"ahler  manifold $X$ the decomposition of complex valued $C^{\infty}$ differential $k$-forms into types induces the Hodge decomposition for the {de} Rham group $H^k(X,\C)$ equipping this group with a pure weight $k$ Hodge structure. For singular or non-compact complex analytic spaces this is no longer true in general. For instance $H^1(\C^*)$ has rank $1$ while it should have even rank if it would carry a weight $1$ Hodge structure.

Cohomology groups of \emph{algebraic} varieties instead carry a canonical \emph{mixed Hodge structure}, i.e. there is a rationally defined increasing weight filtration so that the $k$-th graded pieces carry a weight $k$ Hodge structure. In the {above} example there is only one weight, namely $2$ and $H^1(\C^*)$ is pure of weight $(1,1)$.
In fact, Deligne \cite{Del71,Del74} constructed a good functorial theory for the cohomology of algebraic varieties. 

For a  smooth variety $U$ the weight filtration can be seen on the level of forms as follows. First choose {a so-called \emph{good} compactification, i.e. a smooth (projective) compactification $X$ such that $D=X\setminus U$ is a divisor with normal crossings. De Rham cohomology of $U$ is the cohomology of the full complex of smooth forms on $U$ but  it can also be calculated   using the subcomplex   of  rational forms having at most logarithmic poles along $D$, and the weight filtration is given by the number of logarithmic poles.}     Indeed, $H^k(U)$ carries a mixed Hodge structure with $W_{k-1}H^k(U)=0$ and where $W_kH^k(U)$ is the image of the restriction $H^k(X)\to H^k(U)$.
In the analytic case, a similar assertion holds provided a K\"ahler compactification $X$ of $U$ exists. {All of these assertions are well known consequences of Deligne's theory. } 

{
In the analytic category we work with  manifolds $U$ Zariski-open in some compact K\"ahler space $X$. For these, \emph{good compactifications} exist as in the algebraic case\footnote{More details can be found in \S~\ref{sec3}.}. The weight filtration of the mixed Hodge structure on $H^k(U)$ may (and indeed does) depend on the compactification as shown by the following example.
\begin{exmple*}
$U=\C^*\times \C^*$ can be \emph{analytically} compactified in two ways: one is  $X=\bP^1\times \bP^1$, a second one is the  compactification $Y$ which is  the total space  of the $\bP^1$-bundle   on an elliptic curve $E$ associated to the non-trivial extension of a trivial line bundle by a trivial line bundle. See \cite[Example 4.19]{PS}.  The first has $W_1H^1(U)=0$ while the second has $W_1H^1(U)\simeq H^1(E)$.  Deligne's results  imply that this would not happen if $X$ and $Y$ can be   dominated by a third smooth projective  compactification: indeed, two \emph{birationally} equivalent compactifications  would give the same mixed Hodge structures. 
\end{exmple*}
In the analytic category we are thus led to introduce  the notion of bimeromorphic equivalence: two smooth K\"ahler compactifications $X$ and $Y$ of $U$ are said to be \emph{bimeromorphically equivalent} if they are dominated by a third smooth K\"ahler compactification $Z$ of $U$. If,   moreover, the dominating bimeromorphic morphisms $Z\to X$ and $Z\to Y$ are projective, we say that $X$ and $Y$ are \emph{projective-bimeromorphically equivalent}.  So $U$ always has a good compactification projectively bimeromorphically equivalent  to $X$, but there may be other good compactifications which are not even bimeromorphically  equivalent  to $X$ as  our example  shows.   However,  Deligne's results imply that for any good K\"ahler compactification $Z$ of $U$ we have that $W_kH^k(U)$ is the image of $H^k( Z)$. Hence, in our example, one still has that $W_1H^1(U)$ is the image of $H^1$ of the compactification in both cases. 
 }

{One can generalize the discussion  to cohomology with values in locally constant coefficients where Deligne's theory does not apply. To motivate why one should consider these}, 
look at the Leray  spectral sequence for a  morphism  $f:Y \to X$ between compact K\"ahler spaces; these involve the terms $H^q(X,R^pf_*\underline{\Q}_Y)$.  Assuming  that there is a non-empty Zariski-open subset $U\subset X$ over which $Y$ and $f$ are smooth,   the sheaf $R^pf_*\underline{\Q}_Y|_U$ is indeed locally constant and its  fibers carry a weight $p$ Hodge structure. In fact these can be assembled to give the prototype of what is called a \emph{variation of weight $p$ Hodge structure} (cf. for instance \cite{CSP}). 
So it is  natural to look at $H^k(U,\bV)$ where $\bV$ is a local system.  {The replacement for $H^k(X)$ is intersection cohomology $IH^k(X,\bV)$, and there is an    intrinsic way to relate this  to ordinary cohomology. Indeed,  the adjunction  morphism gives  a  canonical  map $IH^k(X,\bV)\to H^k(U,\bV)$ (see Remark~\ref{WhereItComesFrom} for details).  As is the  case for $R^pf_*\underline{\Q}_Y|_U$, one assumes that $\bV$ carries a  variation of Hodge structure. An extra technical assumption on $\bV$ has to be made which is known to hold\footnote{Because this  system is defined over $\Z$, see \cite[Lemma 4.5]{Sch73}.}  for $R^pf_*\underline{\Q}_Y|_U$: we say that $\bV$ is \emph{quasi-unipotent %local monodromies 
at infinity with respect to $X$} if for some (or any) embedded resolution $X'$ of $(X,X\setminus U)$  the local monodromy operators around the branches of $X' \setminus U$  are   quasi-unipotent.    Indeed, if $\bV$ is quasi-unipotent at infinity 
with respect to  $X$  we have canonical (pure\footnote{its weight is $k+n$ where $n$ is the weight of the variation of Hodge structure  on $\bV$.}, respectively  mixed) Hodge structures on $IH^k(X,\bV)$, respectively $H^k(U,\bV)$. Moreover,  the mixed Hodge structure on $H^k(U,V)$ depends
only the projective bimeromorphic equivalence class of $X$. This will be recalled in \S~\ref{sec3}. See in particular Coroll.~\ref{PureAndMixed}.}

To motivate the statement of the main theorem below, recall Zucker's construction \cite{Zuc79} for $\dim X=1$.
Let $j:U\into X$ be the embedding of $U$ into its compactification. The sheaf $j_*\bV$ is quasi-isomorphic to the complex of holomorphic forms with values in $\bV$ and with $L^2$ growth conditions at the boundary (with respect to the Poincar\'e metric). Forgetting the growth conditions gives a complex which computes the cohomology of $\bV$ on $U$; whence a natural restriction map $L^2H^k(U,\bV) \to H^k(U, \bV)$.  {The source is nothing but another incarnation of $IH^k(X,\bV)$ (Remark~\ref{L2andIC}) and indeed, one of the main results from \cite{Zuc79} states that it has a pure  Hodge structure; moreover,  it  maps  to the lowest weight part of  the (special case of our) mixed Hodge structure on the target.  }
Hence, in this setting, the lowest weight ``comes from the compactification''.  
  
    {The main result of this note concerns a generalization of the lowest weight property containing both Zucker's result and the   constant coefficients case as special cases:
\begin{thm*} \it Assume $U$ is a smooth complex manifold, $j:U\into X$ an analytic-Zariski open inclusion into a compact K\"ahler space  and let $\bV$ be a local system on $U$, {quasi-unipotent at infinity  with respect to $X$ and  carrying a polarizable variation of Hodge structure}.  Then
 \newline
{\rm a)}  {the natural morphism    
 \begin{equation*} 
  IH^k(X, \bV) \to H^k( U, \bV) \eqno{(*)}
 \end{equation*}
 is a morphism of mixed Hodge structures;\\
{\rm b)}   the image of  the map in a)  is exactly the lowest weight part of $H^k( U, \bV)$ and  is the same for  K\"ahler compactifications which are projective-bimero\-morphically  equivalent. }
\end{thm*}
}
{ Let us make some comments on the statement of the theorem and  its proof.}  {Note that $X$ is not assumed to be a \emph{smooth} compactification and that $X\setminus U$ need not be a divisor with normal crossings.  A condition like $X$ being  K\"ahler   is however unavoidable. As to the proof,  a  first point that  needs to be shown is   that the natural map (*)  preserves Hodge and weight filtrations. The second point is that the image of this map, which then lands into the lowest weight part, is exactly the lowest weight part. Finally, since as we have seen, the mixed Hodge structure on $H^k( U, \bV)$ depends on the compactification, one would like to show the lowest weight part depends only on  the bimeromorphic equivalence class, as in the case of constant coefficients.
We can show this only for  projective-bimeromorphically equivalent compactifications. The reason is   that the decomposition theorem in the analytic setting at the moment is only available for projective morphisms (see  \cite[5.3.1]{Sa88}). 
}

{The argument is  not too hard if the complement is a hypersurface. This case is treated in \S~\ref{sec3} {together with the algebraic case}. The general case can be reduced to it by blowing up. This yields  the  mixed Hodge structure on $H^k(U,\bV)$. To see how it relates to  the Hodge structure on $IH^k(X,\bV)$ one needs  a theory of mixed Hodge complexes on analytic spaces developed in \S~\ref{sec4}  and one must study  the direct image of the mixed Hodge complex under  the blowing-down  map. This is done in \S~\ref{sec5} where the final step of the proof is given.}

In order to make the proof as self-contained as possible, we start with a brief summary of the necessary results from the theory of perverse sheaves and mixed Hodge modules.

{We want to make several  side remarks for  the algebraic case.}  The proof we give uses an argument which resembles the one from \cite[Remarks 2.2. i)]{Ha-Sa} used in the $l$-adic situation and for constant coefficients (actually this works as long as the formalism of mixed sheaves \cite{Sa91} is satisfied). 
Note also that  our main theorem in the algebraic case does not follow from the mixed Hodge version of \cite[3.1.4]{Mor} unless $j$ is an affine morphism since the $t$-structure in loc.~cit.\ is {\it not} associated to the mixed complexes of weight $\le k$ in the usual sense, see \cite[3.1.2]{Mor} (and Remark~\ref{Morel} below).

The first named author wants to thank Stefan M\"uller-Stach for asking this question and urging him to write down a proof.

\section{Perverse sheaves} \label{sec1}
 We only give a minimal exposition of the theory of perverse sheaves to explain the properties which will be used below. We shall only be working with the so-called middle perversity which respects Poincar\'e duality. Full details can be found in \cite{B-B-D}.

Let $X$ be a complex analytic space. The category of perverse ``sheaves'' of $\Q$-vector spaces on $X$, denoted by $\perv (X;\Q)$, is an \emph{abelian category}. The fact that it is abelian follows from its very construction as a core with respect to a $t$-structure. While the details of this are not so relevant for what follows, one needs to know that the starting point is formed by the \emph{constructible sheaves of $\Q$-vector spaces} on $X$. By definition these are sheaves of finite dimensional $\Q$-vector spaces which are locally constant on the strata of some analytic stratification of $X$. We assume that the stratification is {\it algebraic} in the algebraic case. The simplest examples of such sheaves are the locally constant sheaves on $X$ itself, or those which are locally constant on some locally Zariski closed subset $Z$ of $X$ but zero elsewhere.

A core is defined with respect to a so-called $t$-structure and in the perverse situation the $t$-structure is defined by certain cohomological conditions, the so called support and co-support conditions. Indeed, instead of starting from complexes of constructible sheaves on $X$ one departs from
\begin{equation} \label{DbC}
\begin{array}{cl}  D_c^{\rm b}(X;\Q):&\text{the derived category of bounded complexes}\\
&  \text{of sheaves of $\Q$-vector spaces on $X$ with} \\
&  \text{constructible cohomology sheaves.}
\end{array}
\end{equation}
By definition a perverse sheaf is such a complex $F$ which obeys the support and co-support conditions:
$$\dim\supp H^{p}(F)\le -p,\quad \dim\supp H^{p}(\D F)\le -p,$$ where $\D F:=\R\underline{\Hom}(F,\D_X)$ is the Verdier dual of $F$ and  $\D_X$ is the dualizing complex. For $X$ smooth and $d$-dimensional{, the case we shall be interested in, $\D_X$ is } just  $\Q_X(d)[2d]$.
The support condition implies that $H^p(F)=0$ for $p>0$ while the co-support condition implies $H^p(F)=0$ for $p<-d$ (where $d=\dim X$): perverse sheaves are complexes ``concentrated in degrees between $-d$ and $0$''.

On a complex manifold a (finite rank) local system of $\Q$-vector spaces $\bV$ can be made perverse by placing it in degree $-d$: the complex $\bV[d]$ is a perverse sheaf. If $X$ is no longer smooth this complex has to be replaced by the so-called intersection complex. Indeed, if $U\subset X$ is a dense Zariski-open subset of $X$ which consists of smooth points and $\bV$ is any (finite rank) local system of $\Q$-vector spaces on $U$ the \emph{intersection complex} $\ic_X(\bV[d])$\footnote{Some people write $\ic_X(\bV)$ instead of $\ic_X(\bV[d])$.} can be constructed as in \cite{B-B-D} (and \ref{ICandPerv} below).
(It is also called the minimal extension.)
By definition, its hypercohomology groups are the intersection cohomology groups:
\begin{equation}\label{eqn:thisisintcoh}
IH^k(X,\bV[d]):= \bH^k(\ic_X(\bV[d])).
\end{equation}

\begin{rmq} Even if $X$ itself is smooth an intersection complex on $X$ need not be of the form $\widetilde\bV[d]$ for some local system $\widetilde\bV$ defined on $X$ because of non-trivial monodromy ``around infinity'' $X\setminus U$.\end{rmq}
The following two results explain the role of these intersection complexes.
\begin{thm}[\protect{\cite[Chap.V, 4]{Bor}}] \label{ICandPerv} Let $X$ be a $d$-dimensional irreducible complex analytic space and let $U$ be a smooth dense Zariski-open subset of $X$ on which there is a local system $\bV$ of finite dimensional $\Q$ vector spaces. The intersection complex $\ic_X(\bV[d])$ is up to an isomorphism in the derived category the unique complex of sheaves of $\Q$-vector spaces on $X$ which is perverse on $X$, which restricts over $U$ to $\bV[d]$ and which has no non-trivial perverse sub or quotient objects supported on $X\setminus U$.
\end{thm}
{
\begin{rmk} \label{WhereItComesFrom} In the situation of  Theorem~\ref{ICandPerv}, let $j:U\into X$ be the inclusion and let $I=\ic_X(\bV[d])$. The adjunction morphism $j^\#: I \to j_*j^* I$ induces a homomorphism
\begin{equation}\label{eqn:WhereItComesFrom}
H^kj^\#: IH^k(X,\bV) \to H^k(U,\bV)
\end{equation}
which will be used to compare intersection and ordinary cohomology.
\end{rmk}
}
\begin{thm}[\cite{B-B-D}] \label{SimplePerverseComplex} If $X$ is compact or algebraic, $\perv(X;\Q)$ is Artinian and Noetherian. Its simple objects are the intersection complexes $F=\ic_Z(\bV[\dim Z])$ supported on an irreducible subspace $Z\subset X$ and where $\bV$ is associated to an irreducible representation of $\pi_1(U)$, $U\subset Z$ the largest open subset of $Z$ over which $F$ is locally
constant.
\end{thm}
We also need \emph{filtered} objects in the abelian category $\perv(X;\Q)$. A priori these are not represented by filtered complexes in the usual sense of the  word, since the morphisms  are in a derived category: they  are ``fractions" $[f]/[s]: K \to L$ where the bracket stands for the corresponding homotopy class, $f:K\to N$ is a morphism of complexes and $N\mapleft{s} L$ is a quasi-isomorphism. However, the category of  sheaves on $X$ with constructible cohomology has enough injectives and replacing $L$ by a complex $L'$ of injective objects, the quasi-isomorphism $s$ becomes invertible up to homotopy and so $[f]/[s]$ can be represented by a true morphism $K  \to L'$.  Next, recall:
\begin{lemma}\label{inj} For any morphism of complexes $v:A\to B$, the morphism in the derived category defined by it can be represented by an injective morphism $A\to B'$ of complexes.\end{lemma}
\proof Take  $B':=\cone(-\id\oplus v:A\to A\oplus B)$. Then  $A$ is a  subcomplex of $B'$ and we get an injective morphism $A\to B'$ which is identified with $v$ by the quasi-isomorphism  $(0,v,\id): B'\to B$. \qed\endproof

\begin{corr} \label{CorInj} Let  $K\in \perv(X;\Q)$. Any  \emph{finite} filtration on $K$ can be represented by a filtered complex in $\perv(X;\Q)$.
\end{corr}
\proof Induction on the length of the filtration, assumed to be an increasing filtration $W$. The above discussion shows that  the morphism $W_i\to W_{i+1}$ in $\perv(X;\Q)$ can be represented by a morphism of complexes to which Lemma.~\ref{inj} can be applied.
\qed\endproof

 \section{Mixed Hodge Modules} \label{sec2}

In this section we put together some properties of mixed Hodge modules which will be used in the sequel.
These properties are proven in \cite{Sa88} and \cite{Sa90}. See also the exposition \cite[Cha. 14]{PS} where mixed Hodge Modules are introduced axiomatically. \medskip

Let $X$ be a complex algebraic variety or a complex analytic space. There exists an abelian category ${\mhm}(X)$, the category of \emph{mixed Hodge modules} on $X$. 
\begin{rmq} Note that for nonproper complex algebraic varieties $X$ we always have ${\mhm}(X)\ne{\mhm}(X^{\rm an})$ because of the difference between algebraic and analytic stratifications. Note also that a mixed Hodge module on an algebraic variety is always \emph{assumed} to be extendable under an open immersion. The last property cannot be well-formulated in the analytic case due to the defect of the Zariski topology on analytic spaces, e.g. Zariski-open immersions are not stable by composition and closed subspaces are not intersections of hypersurfaces Zariski-locally.
\end{rmq}
\begin{proprs} \label{Props}
\begin{lijstb}
\item There is a functor
\begin{equation}
\rat_X:D^{\rm b}\mhm(X)\to D^{\rm b}_c(X;\Q).
\label{eqn:Rat}
\end{equation}
such that $\mhm(X)$ is sent to $\perv(X;\Q)$. One says that $\rat_X M$ is the underlying
rational perverse sheaf of $M$. Moreover, we say that
\[
M\in \mhm(X) \mbox{ \emph{is supported on }} Z \iff \rat_X M \, \mbox{is
supported on } Z.
\]

\item The category of mixed Hodge modules supported on a point is the category
of graded polarizable rational mixed Hodge structures; the functor ``$\rat$'' associates to the mixed
Hodge structure the underlying rational vector space.
\item Each object $M$ in ${\mhm}(X)$ admits a \emph{weight filtration} $W$ such that
\begin{itemize}
\item morphisms preserve the weight filtration strictly;
\item the object $\Gr^W_kM$ is semisimple in $\mhm(X)$;
\item if $X$ is a point the $W$-filtration is the usual weight filtration for the mixed
Hodge structure.
\end{itemize}
Since $\mhm(X)$ is an abelian category, the cohomology groups of any complex of mixed
Hodge modules on $X$ are again mixed Hodge modules on $X$. With this in mind, we say that
for {a} complex $M \in D^{\rm b}\mhm(X)$  the  \emph{weight} satisfies
\[ \mbox{weight}[M]
\left\{ \begin{array}{ll}
\le n, \\
\ge n
\end{array} \right. \iff \Gr_k^WH^i( M)=0 \,\,
\left\{ \begin{array}{ll}
\mbox{for } k>i+n\\
\mbox{for } k<i+n.
\end{array} \right.
\]
We observe that if we consider the weight filtration on the mixed Hodge modules which constitute a complex $M\in D^{\rm b}\mhm(X)$ of mixed Hodge modules we get a filtered complex in this category.
\item\hspace{-0.4em}(i)  For each morphism $f:X\to Y$ between \emph{complex algebraic varieties}, there are induced functors $f_*,f_! :D^{\rm b}\mhm(X)\to D^{\rm b}\mhm(Y)$ and $f^*,f^! :D^{\rm b}\mhm(Y)\to D^{\rm b}\mhm(X)$ which lift  {the} functors $Rf_*$, $f_!$ and $f^{-1}, f^!$ respectively; the latter {two} functors are defined on the level of complexes  {of sheaves on $Y$ (whose cohomology is constructible.)}

\noindent D)(ii) In the \emph{analytic case}  {(i)}  is no longer necessarily true but we have:
\\
 | for   $f:X\to Y $  projective or if $X$ is compact K\"ahler and $Y=\text{\rm pt}$, there are \emph{cohomological} functors $H^if_*=H^if_! :\mhm(X)\to\mhm(Y)$ which lift  the perverse cohomological functor $^pR^if_*={}^pR^if_!$;
\\
| for any $f$ there are cohomological functors $H^if^*,H^if^! :\mhm(Y)\to\mhm(X)$ which lift $^pH^if^{-1}$, $^pH^if^!$ respectively;
\item The functors $f^*,f_!$ do not increase weights in the sense that if $M$
has weights $\le n$, the same is true for $f^*M$ and $f_!M$.

\item The functors $f_*,f^!$ do not decrease weights in the sense that if $M$ has
weights $\ge n$, the same is true for $f_*M$ and $f^!M$.

\item If $f$ is proper, $f_*$ preserves weights, i.e. $f_*$ neither increases nor decreases weights.
\end{lijstb}
\end{proprs}
\begin{rmks}\label{rmks2}
1)  Despite the fact that the functors $f^*$ etc. do not exist in the analytic setting, properties E), F), G) still have a meaning as in \cite[2.26]{Sa90} since the weight is defined in terms of  cohomology  only.
\newline
2)   Since in the analytic setting Zariski-open immersions are not stable by composition $H^if_*M$, $H^if_!M$ do not necessarily exist for analytic morphisms $f$. This explains why in the analytic case D) not all morphisms are allowed.
\newline
3) The reader may interpret   the K\"ahler condition on $X$ in Property D) as  the existence of  a projective morphism $g$ from a K\"ahler \emph{manifold} $X'$ onto $X$. Indeed, the construction of $H^if_*M$ for $f:X\to $pt , where $M$ is a pure Hodge module,  is reduced to the assertion for $X'$: use the decomposition theorem for $g$ applied to a pure Hodge module on $X'$ which is a subquotient of the pullback of $M$ by $g$. Then it follows from \cite{CKS}, \cite{KK86}, \cite{KK87}, \cite{KK89}. For the mixed case we can use the weight spectral sequence. 
\newline
4)  It is still unclear whether $H^if_*M$  exists for proper K\"ahler morphisms $f$ unless $M$ is constant, see \cite{Sa90b}. If the reader prefers, he may assume that the polarizable Hodge modules in this paper are direct factors of the cohomological direct images of the constant sheaf by smooth K\"ahler morphisms so that   {the existence of $H^if_*M$}  follows from the decomposition theorem for the direct image of the constant sheaf by proper K\"ahler morphisms \cite{Sa90b}.
\end{rmks}

The above properties readily imply various basic properties of mixed Hodge modules. For example,
if $M$ is a complex of mixed Hodge modules on $X$ its cohomology $H^q M $ is a mixed Hodge module on $X$.
Properties B) and D) imply:
\begin{lemma} \label{HyperCohForMHS} Let $a_X:X\to \mbox{\rm pt}$ be the constant map to the point. Assume $X$ is algebraic or compact K\"ahler. Then for any complex $M$ of mixed Hodge modules on $X$
\begin{equation} \label{eqn:HyperCohForMHS}
\bH^p(X,M) := H^p((a_X)_* M)
\end{equation}
is a mixed Hodge structure.
\end{lemma}
For the proof of the main theorem one needs the following two technical constructions. The first is the adjunction construction:
\begin{constr}
Consider a morphism $f:X\to Y$ of {\it algebraic varieties} and a mixed Hodge module $M$ on $Y$. The adjunction morphism $f^\#: M \to f_* f^* M$ is a morphism of complexes of mixed Hodge modules. For any bounded complex $K$ of mixed Hodge modules on $X$, the identity $a_X=a_Y\comp f$ induces a canonical identification $\bH^n(Y, f_* K)= \bH^n(X, K)$. In particular this holds for $K=f^*M$. Adjunction thus induces a morphism of mixed Hodge structures
\begin{equation} \label{eqn:AdjMorph} H^kf^\#: \bH^k( Y, M) \to \bH^k(X, f^* M).
\end{equation}
\end{constr}
In the \emph{analytic case}  this construction remains valid   for an open immersion $j$  whose complement is a hypersurface (defined locally by a function $g$). Indeed, then $j_*j^*M$ is a mixed Hodge module whose underlying $D$-module  comes from  localization by $g$  {applied to the underlying $D$-module of $M$}.  More generally, consider  the complement $U$ of an intersection $Z$ of global hypersurfaces. Then $j_*j^*M$ is a complex of mixed Hodge modules  due to a  second construction:
\begin{constr}[\protect{\cite[2.19, 2.20]{Sa90}}] \label{construction}  Let $g_i$, $i=1,\dots,r$ be holomorphic functions on $Y$, let $Z=\bigcap_{i=1}^r g_i^{-1}(0)_{\rm red}$ and $U= Y \setminus Z$. We set  $Y_i= Y \setminus  g_i^{-1}(0)$ and  for  $I\subset\set{1,\dots,r}$  we set $Y_I =\bigcap _{i\in I} Y_i $. Let 
$i: Z\into Y $, $j: U\into Y$ and  $ j_I:Y_I\into Y $ be the natural inclusions.  Let $M$ be a mixed Hodge module on $Y$;  then  also $j^*M$, the restriction of $M$ to $U$, is a mixed Hodge module on $U$ and    there are quasi-isomorphisms in the category $D^{\rm b} \mhm(Y)$
\begin{eqnarray*}
i_*i^! M & \mapright{\sim} &[ \cdots 0 \to M \to B_1\to B_2   \cdots B_r\to 0]  ,\quad B_{k}= \bigoplus_{|I|=k} (j_I)_* j^*_I M \\
j_*j^*M& \mapright{\sim} &  [\cdots 0  \to B_1\to B_2 \to B_3\cdots\to B_r\to 0 ], \quad \text{ $B_{k}$  in degree $ k-1$. }  
\end{eqnarray*}
\end{constr}
The above construction  leads to
\begin{lemma}[\protect{ \cite[(4.4.1)]{Sa90}}] \label{dtriangle} Let $i: Z\subset Y$ be a closed immersion and $j: U=Y\setminus Z\into Y$ be the inclusion of the complement. Assume $Y,Z$ are algebraic or, alternatively, that $Z$ is an intersection of global hypersurfaces of $Y$. Let $M$ be a  
mixed Hodge module on $Y$. There is a distinguished triangle \footnote{We shall write   triangles  also as $M' \to M \to M''\to [1] $. }
\begin{equation} \label{eqn:dtriangle}
\hspace{6em} \begin{diagram}
 i_*i^! M & \mapright{\quad \quad} & M \\
 \vspace{3ex} \\
 & j_*j^{*}M \\
 \end{diagram}
 \setlength{\unitlength}{1.1mm}
 \begin{picture}(42,4)(0,0)
 \arrow(-20,-4)\dir(-3,4)\long{5}
 \arrow(-5,5)\dir(-3,-4)\long{5}
 \put(-7,0){$\scriptstyle \alpha$}
 \put(-22,0){$\scriptstyle [1]$}
 \end{picture}\hspace{-6em}
 \end{equation}
 in the bounded derived category of mixed Hodge modules lifting the analogous triangle for complexes with constructible cohomology sheaves. The morphism $\alpha$ induces
 the adjunction morphism $H^kj^\#: \bH^k(Y,M) \to \bH^k(U, j^*M)$ for $j$ (see \eqref{eqn:AdjMorph}).
\end{lemma}
\proof
In the algebraic setting the  local constructions~\ref{construction}  for a suitable affine cover {patch}  together to give  globally defined quasi-isomorphisms for  $i_*i^!M$ and $j_*j^*M$. The local construction shows the existence of the distinguished triangle. See  \cite[4.4.1]{Sa90} for details.
 
The same argument applies in the analytic case under the assumption that  $Z$  is a global complete intersection. See  the proof of \cite[2.19]{Sa90}.   \qed\endproof
\begin{rmk} The reader may wonder what happens in the general setting of analytic spaces. The problem is that Construction~\ref{construction} can not be globalized to complexes of mixed Hodge modules. However, the \emph{cohomology sheaves} of the complexes do make sense globally and are indeed mixed Hodge modules. Hence also the long exact sequence in cohomology associated to \eqref{eqn:dtriangle} exists in the category of mixed Hodge modules. 
\end{rmk}

\section{Polarizable variations of Hodge structure and the main  theorem} \label{sec3}

In this section $X$ is an \emph{irreducible} compact K\"ahler analytic space of dimension $d$.
Let $j:U\into X$  {be} the inclusion of a dense Zariski-open subset for which we make the crucial assumption that  \textbf{ \textit{$U$ is smooth.}}

{We shall not review here the definitions and properties of polarizable Hodge modules.
For our purposes we only need  the following  basic result linking variations  of Hodge structures  and polarizable Hodge modules} \cite[Th.~5.4.3]{Sa88}:
\begin{thm} \label{PolHS} Suppose that $\bV$ is a polarizable variation  {of Hodge structure} on $U$ of weight $n$. {If  $U$ is smooth}, there is a polarizable Hodge module $V^{\rm Hdg}$ of weight $n+d$ on $U$ whose underlying perverse component is $\bV[d]$.
\end{thm}
{There is however, one important aside to make at this point.
In the algebraic setting mixed  Hodge modules are assumed to extend under  open immersions}{, but this ceases to hold in the analytic category.  Instead, one replaces it  by the following  condition on   the underlying local system. }
\begin{dfn} { Let $\bV$ be a local system on $U$. We say that $\bV$  is \emph{quasi-unipotent at infinity  with respect to $X$},  if  for some (or any) choice of an embedded resolution $X'$ of $(X,X\setminus U)$  
 the   local monodromy operators of $\bV$ around $X'\setminus U$ are quasi-unipotent.
} 
\end{dfn}
\begin{rmq}{By   \cite{Kas} this  property is independent of the choice of $X'$ and depends only on  the bimeromorphic equivalence class of  $X$.  By definition $X'$  is a good compactification of $U$  together with a  bimeromorphic map $f:X'\to X$; these exist: by  blowing up $X$ in suitable ideals on can even assume that $f$ is projective.  So one can test quasi-unipotency on such $X'$.
}
\end{rmq}
For a polarizable Hodge module  this {notion} leads to   \emph{pure} Hodge modules (see  \cite{Sa90}). 
{In both settings  (polarizable algebraic Hodge modules and pure mixed Hodge modules in the analytic case)   one obtains} semi-simple categories: this is implied by the polarizability condition, see also \ref{Props}.C). {Both categories satisfy} moreover the strict support condition:
\begin{propr} A polarizable weight $n$ Hodge module $M$ is a direct sum of polarizable weight $n$ Hodge modules $M_Z$ which have strict support\footnote{$M$ is said to have \emph{strict support $Z$} if it is supported on $Z$ but no quotient or sub object of $M$ has support on a proper subvariety of $Z$.}  $Z$ where $Z$ are irreducible subvarieties of $X$, and the same assertion holds for pure Hodge modules.
\end{propr}
By \cite[3.20, 3.21]{Sa90} one has:
\begin{thm} \label{PolHM}  Assume that $U$ is smooth and   that $\bV$ is a quasi-unipotent  at infinity  with respect to $X$  and underlies a polarized variation of Hodge structures of weight $n$ on $U$. Then there is a unique pure Hodge module $V^{\rm Hdg}_X$ of weight $n+d$ on $X$ having strict support in $X$ and which restricts over $U$ to $V^{\rm Hdg}$.
\end{thm}
\begin{rmq} Note that this checks with the assertion in Theorem~\ref{ICandPerv} which holds for the rational component of the mixed Hodge modules.
 {
 More precisely: the intersection complex $\ic_X(\bV[d])$ is the rational component of  $V^{\rm Hdg}_X$. This remark is crucial for the proof of the next corollary. }
\end{rmq}
  \begin{corr} \label{PureAndMixed} {  {\rm 1)} There exists  a mixed Hodge structure on $H^k(U,\bV)$. It depends only on the projective bimeromorphic equivalence class of $X$; \\
{\rm 2)} $IH^k(X,\bV)$ carries a pure Hodge structure of weight $k+n$. }
\end{corr}
 \proof
{1) Replacing $X$ by a suitable blow up, we may  assume that $X$ is a good compactification of $U$. By construction~\ref{construction} $j_*j^*V_X^{\rm Hdg}=j_*V^{\rm Hdg}$ then is a  mixed Hodge module} on $X$ and by Lemma~\ref{HyperCohForMHS} the cohomology group $H^k(U,\bV)$ carries the  mixed Hodge structure $\bH^{k-d}(j_*V^{\rm Hdg})$.  If there are two good compactifications $X_1,X_2$ with a projective morphism $\pi:X_1\to X_2$ inducing
an isomorphism over $U$, then with $j_k:U\into X_k$, $k=1,2$ the embeddings, we have an isomorphism of mixed Hodge modules $R\pi_*Rj_{1*}V_{X_1}^{\rm Hdg}=Rj_{2*}V_{X_2}^{\rm Hdg}$ by the uniqueness of $Rj_*$, see e.g. \cite[2.11]{Sa90}. 
\\
2)  Since $V_X^{\rm Hdg}$ is a pure Hodge module, by the previous Remark and Lemma~\ref{HyperCohForMHS} $ H^k(X,\ic_X(\bV)=IH^k(X,\bV)$ carries the mixed Hodge structure  $H^{k-d}(X,V_X^{\rm Hdg})$, which by Properties~\ref{Props} G) is pure of weight $k+n$.
 \qed 
 \endproof
\begin{rmk} \label{L2andIC}
Suppose that in addition $X$ is smooth and $X-U$ is a divisor with normal crossings. Then, by \cite[Theorem 1.5]{CKS}, \cite{KK86}, \cite{KK87}, \cite{KK89} $IH^k(X,\bV)$ can be identified with $L^2H^k(U,\bV)$ provided one measures integrability with respect to the Poincar\'e metric around infinity (one is in the normal crossing situation, so locally around infinity {one has}  a product of disks and punctured disks).
Summarizing:
\[
\bH^k(\ic_X(\bV)) = IH^k(X, \bV)= L^2H^{k}(U ,\bV)
\]
has a pure Hodge structure of weight $k+n$.
\end{rmk}
 
Next one wants to relate intersection cohomology and ordinary cohomology. This is the content of the main theorem:
\begin{thm} \label{mainresult} {Assume   that $U$ is smooth, that $\bV$ is a quasi-unipotent at infinity  with respect to $X$ and that it  carries a polarized variation of Hodge structure of weight $n$}.
 {Then
 \newline
{\rm a)} The natural  morphism 
 \[
H^kj^\#: IH^k(X, \bV) \to H^k( U, \bV) 
\] 
(see  \eqref{eqn:WhereItComesFrom}) is a morphism of mixed Hodge structures;}
 \\
{\rm b)}   the image of $H^kj^\#$   is exactly the lowest weight part of $H^k( U, \bV)$ {and this image  is the same for  K\"ahler compactifications which are projective-bimero\-morphically  equivalent to $X$.}
\end{thm}
{
\noindent \textit{Proof in the algebraic case.} }
Let $i:Z= X\setminus U \into X$ be the inclusion.     Set  $M=V^{\rm Hdg}_X$, $M'=j_*V^{\rm Hdg}= j_*j^* V^{\rm Hdg}_X$ and   $M''=  i_*i^! V^{\rm Hdg}_X$.  
{Formula \eqref{eqn:AdjMorph} for the inclusion $j:U\into X$ and the mixed Hodge module $M:=V_X^{\rm Hdg}$ shows that \eqref{eqn:WhereItComesFrom} is indeed a morphism of mixed Hodge structures.} 

Form  the distinguished triangle \eqref{eqn:dtriangle}. Portion of its associated long exact sequence in hypercohomology reads
\begin{equation}\label{eqn:crucial}
\hspace{-0.5em}\begin{matrix}
\cdots\to  IH^k(X,\bV)  & \mapright{H^kj^\#} & H^k(U,\bV)&& \\ 
\Vline{}{}{2ex}  && \Vline{}{}{2ex} &&\\
 \bH^{k-d}(X,M)  &\mapright{\phantom{H^kj^\#}} &  \bH^{k-d}(X,M') & \to & \bH^{k-d+1}(X, M'') \to \dots
\end{matrix}
\end{equation}
By Theorem~\ref{PolHM} $M=V^{\rm Hdg}_X$ is pure of weight $n+d$, and so by Property~\ref{Props}.F the complex $i^!V^{\rm Hdg}_X$ has weight $\ge n+d$. By Property~\ref{Props}.G this also holds for the complex $M''=i_*i^!V^{\rm Hdg}_X$. Applying once more Property~\ref{Props}.G to the functor $(a_X)_*$ one sees that $\bH^{k-d+1}(X, M'') $  has weights $\ge k+n+1$ and hence the image of the map \eqref{eqn:WhereItComesFrom} is exactly the weight $(k+n)$-part of $H^k( U, \bV)$. 

{In the algebraic category  the last assertion of b) can be replaced by a  stronger assertion: we may assume that  two compactifications are related by a proper algebraic morphism to which the decomposition theorem can be applied. Instead of giving full details here we refer to the proof in the analytic setting which is given at the end of  \S~\ref{sec5} and which is similar in spirit.}

{\noindent \textit{Strategy of the proof in the analytic setting.}}
If   $Z$ is a  hypersurface,  the same proof works in view of Lemma~\ref{dtriangle}. {This Lemma also shows that  \eqref{eqn:WhereItComesFrom}  is a morphism of mixed Hodge structures in this case. }

In the general situation one has to perform a suitable blow-up $\pi:X'\to X$   which is the identity in $U$ and such that $Z'=X'\setminus U$ is a divisor. Now we would like to apply the functor $\pi_*$. The problem is that this functor does not exist in  the derived categories of mixed Hodge modules. So we have to find a substitute for this which still preserves  enough of the Properties ~\ref{Props} so that we can complete the proof as in the algebraic case. It turns out that the correct category to use is the one of mixed Hodge complexes. See \S~\ref{sec4}.    {In  \S~\ref{sec5} we complete the proof in the analytic case.}  

\begin{rmk} \label{finalrmk} In the  algebraic setting  the following claim is easily shown to imply the main  result as well and can be seen as a refinement of it.
 \begin{claim} Suppose $Z$ is a locally principal divisor or $j$ is an affine morphism. Then the adjunction morphism $j^\#: V^{\rm Hdg}_X \to j_*j^*V^{\rm Hdg}_X$ is injective and identifies $V_X^{\rm Hdg}$ with the lowest weight part of $ j_*j^*V^{\rm Hdg}_X=j_* V^{\rm Hdg}$
\end{claim}
Indeed, the extra hypothesis on $j$ implies  
(see Construction~\ref{construction} and \cite[2.11]{Sa90}) 
 that $j_* V^{\rm Hdg}$  is a mixed Hodge module (not just a \emph{complex} of mixed Hodge modules) and the main theorem then follows easily from the Claim.  {The latter  follows  from the  long exact sequence $0\to H^0i_*i^!V^{\rm Hdg}_X\to V^{\rm Hdg}_X\to j_* V^{\rm Hdg}\to H^1i_*i^!V^{\rm Hdg}_X$ using that  the strict support condition implies that $H^0i_*i^!V^{\rm Hdg}_X=0$.}

{The above claim can alternatively be shown using adjunction. This was how the first named author originally proved the main result\footnote{See {\tt http://arxiv.org/abs/0708.0130v2} }.  Here is the argument. It  suffices  to show that the  lowest weight part $W_{d+n}M$ of $M=j_*V^{\rm Hdg}$ has  no quotient or  sub object supported on $D=X\setminus U$.  It is  a  pure weight mixed Hodge module, and hence, by Property~\ref{Props}.C  a semi-simple object in the category of mixed Hodge modules. By construction, it restricts to $V^{\rm Hdg}$ on $U$.  By semi-simplicity a quotient object is also a sub object and hence it suffices to show that there are  no    mixed Hodge  modules   $N$   of pure weight $n+d$ supported on $D$ for which $\Hom_A (N,  W_{n+d} M)=0$ in the abelian   category    $A$ of mixed Hodge modules. Functoriality of the weight filtration implies 
$\Hom_A (N,  W_{n+d} M)=\Hom_A(N,M)$.   Let $D(A)$ be the derived category of bounded complexes in $A$. Since the natural map $\Hom_A(N,M)\to \Hom_{D(A)}(N,M)$ is a bijection (see \cite[p. 293]{Verd77}) it is enough to show that $\Hom_{D(A)}(N,M)=0$. In the derived category one  can use the adjunction for $(j^*,j_*)$ yielding  $\Hom_{D(A)}(N,j_*V^{\rm Hdg})=\Hom_{D(A)}(j^*N,V^{\rm Hdg})=0$ since $j^*N=0$.}
\end{rmk}

\begin{rmk}\label{Morel} It would not be difficult to construct a mixed Hodge version of \cite[3.1.4]{Mor}. However, this would not immediately imply our main theorem unless $j$ is an affine morphism. Indeed, the $t$-structure in loc.~cit.\ is defined by the condition that $^pH^iK$ has weight $\le k$ and not $\le i+k$ as in the case of mixed Hodge complexes of weight $\le k$, see \cite[3.1.2]{Mor}. It does not seem that there exists a  $t$-structure associated to mixed complexes of weight $\le k$ since the weight filtration is not strict and the weight spectral sequence does not degenerate at $E_1$ (see also Section~\ref{sec5} on the proof of Theorem~\ref{mainresult} in the analytic case where mixed Hodge complexes in the Hodge setting are used).\end{rmk}

\section{Mixed Hodge complexes on analytic spaces} \label{sec4}

For the proof of Theorem~\ref{mainresult} in the analytic case we need a theory {of} mixed Hodge complexes on analytic spaces which   refines Deligne's theory  \cite{Del71} of  cohomological mixed Hodge complexes. We present  it here in a rather simplified manner which  has the  defect that the mapping cones are not well-defined. However,  this does not cause a problem for the proof of Theorem~\ref{mainresult} since all we need is the existence of the long exact sequence \eqref{eqn:longseq}. See \cite{Sa00} for a  more elaborate formulation taking care of the problem with the cones.
\begin{nota}  
|  $MFW(D_X)$:  the category of filtered $D_X$-modules $(M,F)$ with a finite filtration $W$. For singular $X$ this can be defined by using closed embeddings of open subsets of $X$ into complex manifolds, see \cite[2.1.20]{Sa88}. 
\newline
| $D^{\rm b}_hFW(D_X)$:  the derived category of bounded complexes $(M,F,W)$ such that 
1)  the sheaves  $\bigoplus_pH^iF_p\Gr^W_kM$ are coherent over the sheaf $\bigoplus_pF_pD_X$ 
and 2) the sheaves $H^i\Gr^W_kM$ are holonomic $D_X$-modules.
\newline
 | $D^{\rm b}_cW(X,\Q)$:  the derived category of of bounded filtered complexes $(K,W)$ such that $W$ is finite and $\Gr^W_kK\in D^{\rm b}_c(X,\Q)$ for any $k$: we define $D^{\rm b}_cW(X,\C)$ similarly. 
\newline 
|  $D^{\rm b}_hFW(D_X,\Q)$:  the ``fibre'' product of $D^{\rm b}_hFW(D_X)$ and $D^{\rm b}_cW(X,\Q)$ over $D^{\rm b}_cW(X,\C)$  where  the  functor  $\DR:D^{\rm b}_hFW(D_X)\to D^{\rm b}_cW(X,\C)$ induced by the {de} Rham functor is used to glue the two categories. More precisely,  its objects are triples
$$\MM=((M,F,W),(K,W),\alpha)$$
where $(M,F,W)\in D^{\rm b}_hFW(D_X)$, $(K,W)\in D^{\rm b}_cW(X,\Q)$ and
$$\alpha:\DR(M,W)\cong(K,W)\otimes_{\Q}\C\quad\hbox{in}\,\,D^{\rm b}_cW(X,\C)$$
and morphisms in the category are pairs of morphisms of $D^{\rm b}_hFW(D_X)$ and $D^{\rm b}_cW(X,\Q)$ compatible with $\alpha$. 
 Forgetting the filtration $W$ we can define $D^{\rm b}_hF(D_X)$, $D^{\rm b}_c(X,\Q)$ and $D^{\rm b}_hF(D_X,\Q)$ similarly.  
\newline
$$\text{\hspace{-3em}|  \hspace{3em}}  \Gr^W_k\MM=(\Gr^W_k(M,F),\Gr^W_kK,\Gr^W_k\alpha)\in D^{\rm b}_hF(D_X,\Q).$$
\end{nota}
\begin{dfn}\label{dfn4}
1) The category of mixed Hodge complexes $\mhc(X)$  is the full subcategory of $D^{\rm b}_hFW(D_X,\Q)$ consisting of $\MM=((M,F,W),(K,W),\alpha)$ satisfying the following conditions for  $\Gr^W_k\MM$ for any $k,i$:
\newline
(i) The $\Gr_k^W(M,F)$ are strict and we have a decomposition
\begin{equation}\label{eqn:dec}\Gr^W_k\MM\cong\bigoplus_j(H^j\Gr^W_k\MM)[-j].
\end{equation}
(ii) The $H^i\Gr^W_k\MM$ are polarizable Hodge modules of weight $k+i$.

\noindent 2) Let $\mhw(X)$ denote the category of \emph{weakly mixed Hodge modules}, i.e. its objects have a weight filtration $W$ fo{r} which the gradeds $\Gr^W_k$ are polarizable Hodge modules of weight $k$, but there is no condition on the extension between the graded pieces.

\noindent 3) We say that $\MM\overset{u}\to\MM'\overset{v}\to\MM''\overset{w}\to\MM[1]$ is a \emph{weakly distinguished triangle} in $\mhc(X)$ if $u,v,w$ are morphisms of $\mhc(X)$ and its underlying triangle of complexes of sheaves of $\Q$-vector spaces is distinguished. Here the weight filtration $W$ on $\MM[1]$ is shifted by 1 so that $\MM[1]$ is a mixed Hodge complex.
\end{dfn}
 \begin{rmq}  
In the case $X=$pt , we do not have to assume the decomposition~\eqref{eqn:dec} in condition (i) of Definition~\ref{dfn4},1). One reason is that this is only needed to prove the stability by the direct image under a morphism from $X$. Another reason is that this decomposition actually follows from the other conditions in this case since the category of vector spaces over a field is semisimple.
\end{rmq}
We have by \cite[5.1.14]{Sa88}
\begin{prop}\label{abel} The category $\mhw(X)$ is an abelian category whose morphisms are strictly compatible with $(F,W)$.
\end{prop}

For a mixed Hodge complex $\MM$, set
$$H^i\MM=(H^i(M,F),{}^pH^i(K),{}^pH^i\alpha).$$
We put a weight filtration on it by letting $W_k$ be  the image of $H^iW_{k-i}\MM$ (or, equivalently, the one  induced by the filtration $\Dec W$   for the underlying $D$-module (cf. Proposition~\ref{specseq} below). This shift of the filtration $W$ comes from condition (ii) in the above definition of $\mhc(X)$.

Using \cite[1.3.6 and 5.1.11]{Sa88}, etc. we have
\begin{prop}\label{specseq} With the weight filtration $W$ defined above, the $H^i\MM$ are weakly mixed Hodge modules.
There is a  weight spectral sequence in the abelian category of weakly mixed Hodge modules $\mhw(X)$
\begin{equation}\label{eqn:specseq}E_1^{p,q}=H^{p+q}\Gr^W_{-p}\MM\Rightarrow H^{p+q}\MM,\end{equation} which degenerates at $E_2$, and whose abutting filtration on $H^{p+q}\MM$ coincides with the weight filtration of weakly mixed Hodge modules shifted by $p+q$ as above, i.e.
\begin{equation}\label{eqn:infty}
E_\infty^{p,q}= \Gr^W_{q} H^{p+q}\MM
\end{equation}
Moreover, $(M,F,\Dec W)$ is bistrict, and the weight filtration on $H^{p+q}M$ is induced by $\Dec W$ where $M$ is the underlying $D$-module of $\MM$ and
$$(\Dec W)_kM^i:=\Ker(d:W_{k-i}M^i\to\Gr^W_{k-i}M^{i+1}).$$\end{prop}
Combining this with Proposition~\ref{abel} we get
\begin{prop}\label{longseq} A weakly distinguished triangle as in Definition~\ref{dfn4}, 3) induces a long exact sequence in the abelian category $\mhw(X)$
\begin{equation}\label{eqn:longseq} \to H^i\MM\overset{u}\to H^i\MM'\overset{v}\to H^i\MM''\overset{w}\to H^{i+1}\MM\to.\end{equation}
\end{prop}
For a morphism of mixed Hodge complexes $u:\MM\to\MM'$, there is a {\it mapping cone} $\MM'':=\cone(u:\MM\to\MM')$ in the usual way. Here the weight filtration $W$ on $\MM[1]$ is shifted by 1 so that $\Gr^Wu$ in the graded pieces of the differential of $\MM''$ vanishes and hence conditions (i) and (ii) above are satisfied. However, $\MM''$ is {\it not unique} up to a non-canonical isomorphism because of a problem of homotopy. So we cannot get a triangulated category although there is a weakly distinguished triangle $\MM\to\MM'\to\MM''\to[1]$ 
which  by Proposition~\ref{longseq} induces the long exact sequence \eqref{eqn:longseq} in the category $\mhw(X)$.

Since the weight filtration on the perverse component of a weakly mixed Hodge module can be represented by an honest filtered complex (Cor.~\ref{CorInj}) we have:
\begin{prop}\label{fun} Considering a  weakly mixed Hodge module  as a mixed Hodge complex concentrated in  degree $0$ we get  a functor
$$\iota_X:\mhw(X)\to\mhc(X).$$
\end{prop}
{Let $f:X\to Y$ be a projective morphism, and let $M$ be a polarizable Hodge module. The object $f_*(\iota_X(M))$ belongs to $D^{\rm b}_hF(D_Y,\Q)$. The decomposition theorem \cite[5.3.1]{Sa88} can be applied to $M$ and applying $\iota_Y$ to the resulting Hodge modules yields elements in $D^{\rm b}_hF(D_Y,\Q)$. The uniqueness of the decomposition \cite{Del94} then implies:}
\begin{thm}\label{decomp} Let $f:X\to Y$ be a projective morphism, and $\MM$ be the image of a polarizable Hodge module by $\iota_X$. Then we have a decomposition
$$f_*\MM\cong\bigoplus_i(H^if_*\MM)[-i]\quad\text{in}\,\,\,D^{\rm b}_hF(D_Y,\Q).$$
\end{thm}
Combining this with Properties~\ref{Props}.D) (ii) we get
\begin{corr}\label{stable} Mixed Hodge complexes and weakly distinguished triangles are stable by the direct image under $f:X\to Y$ if $f$ is projective or if $X$ is compact K\"ahler and $Y=pt$.\end{corr}
\begin{rmq}     Note that the stability by direct images  asserted in Corollary~\ref{stable} does not follow from Theorem~\ref{decomp} if we replace $k+i$ by $k$ in condition (ii) in the above definition of $\mhc(X)$. (This causes the shift of the filtration $W$ in Proposition~\ref{specseq} below.)
\end{rmq}

\section{Proof of Theorem~\ref{mainresult} in the analytic case} \label{sec5}

Let $\pi:X'\to X$ be a bimeromorphic projective morphism inducing the identity over $U$ and such that $X'-U$ is a hypersurface (defined locally by a function). Let $j':U\to X'$ denote the inclusion. Then $\R j'_*\bV[d]$ is a perverse sheaf, and underlies a mixed Hodge module $j'_*V^{\rm Hdg}$, see \cite[2.17]{Sa90}. By Proposition~\ref{fun} this gives a mixed Hodge complex concentrated in degree $0$
\begin{equation} \label{eqn:CrucialMHC}
\MM'=((M',F,W),(K',W),\alpha):=\iota_{X'}(j'_*V^{\rm Hdg}),
\end{equation}
such that $K'=\R j'_*\bV[d]$ and $\MM'|_U$ is identified with $V^{\rm Hdg}$.
We denote the direct image of $\MM'$ by
$$\MM=((M,F,W),(K,W),\alpha):=\pi_*\MM'=(\pi_*(M',F,W),\pi_*(K',W),\pi_*\alpha).$$
By Corollary~\ref{stable} this is a mixed Hodge complex  since $\pi$ is projective.

\begin{prop} \label{keyprop} We have $\Gr^W_{d+n}H^0\MM=\iota_X(V_X^{\rm Hdg})$, and $\Gr^W_kH^i\MM=0$ if $k=d+n+i,i\ne0$ or if $k<d+n+i$.\end{prop}
\proof It suffices to show the assertion for the underlying complex of $D$-modules $M$ since the condition on strict support in Theorem~\ref{PolHM} is detected 
 by its underlying $D$-module. 
 Moreover we may restrict to a sufficiently small open subset $Y$ of $X$ 
enabling us to  apply Construction~\ref{construction}.

So let $g_1,\dots,g_r$ be functions on $Y$ such that $Z\cap Y=\bigcap_ig_i^{-1}(0)$. Set $Y_i=Y-g_i^{-1}(0)$.  Abusing notation,    let $i:Y\cap Z\to Y$, $j:Y-Z\to Y$ denote the inclusions.  
By Lemma \ref{dtriangle} there is a  distinguished triangle
$$i_*i^!(V_X^{\rm Hdg}|_Y)\to V_X^{\rm Hdg}|_Y\to j_*j^*(V_X^{\rm Hdg}|_Y)\to [1],$$
inducing a long exact sequence of cohomology. 
\begin{claim} The underlying  bifiltered $D$-modules of  $\iota_Y(H^ij_*j^*(V_X^{\rm Hdg}|_Y))$ and $H^i\MM|_Y$ are isomorphic to each other.
\end{claim}
Suppose that the Claim has been shown. Then the same argument as in the proof of Theorem~\ref{mainresult} in the algebraic case {proves the result of the Proposition}. Indeed, we have the exact sequence
\[
H^i(V_X^{\rm Hdg}|_Y)\to H^ij_*j^*(V_X^{\rm Hdg}|_Y)\to H^{i+1} i_*i^!(V_X^{\rm Hdg}|_Y),
\]
and $H^{i+1}i_*i^!(V_X^{\rm Hdg}|_Y)$ has weights $\ge d+n+i+1$ by Properties~\ref{Props}.F and G. This gives the assertion for $i=0$ since  $V_X^{\rm Hdg}|Y$ is pure of weight $d+n$. For $i\ne 0$ we have $H^i(V_X^{\rm Hdg}|Y)=0$ and hence the last morphism of the exact sequence is injective so that the assertion follows.
\newline
\noindent\textit{Proof of the Claim.}
Let $Y'=\pi^{-1}(Y)$, $Y'_i=\pi^{-1}(Y_i)$, and $g'_i=\pi^*g_i$.  
By Construction~\ref{construction}  
the associated \v Cech complex gives a resolution of $j'_*V^{\rm Hdg}$.
The  components of this \v Cech complex 
are direct sums of $(j'_I)_*(V^{\rm Hdg}|Y'_I)$ where $Y'_I=\bigcap_{i\in I}Y'_i$ with the inclusion $j'_I:Y'_I\to Y$. By the uniqueness of the open direct image in \cite[2.11]{Sa90} we have moreover
$$\pi_*(j'_I)_*(V^{\rm Hdg}|Y'_I)=(j_I)_*(V^{\rm Hdg}|Y_I),$$
where $j_I:Y_I:=\bigcap_{i\in I}Y_i\to Y$. So we get the desired isomorphism (using the filtration $\Dec W$  from  Proposition~\ref{specseq}), and Proposition~\ref{keyprop} follows.\qed\endproof

\medskip
We return to the proof of Theorem~\ref{mainresult} in the analytic case. Applying Proposition~\ref{keyprop} to $\MM'$, we get
$$\Gr^W_{d+n}\MM'=\iota_{X'}(\Gr^W_{d+n}j'_*V^{\rm Hdg})=\iota_{X'}(V_{X'}^{\rm Hdg}).$$
This implies that we get a morphism $u':  \iota_{X'}(V_{X'}^{\rm Hdg}) \to \MM'$ to which we apply $\pi_*$. The decomposition Theorem  \ref{decomp} together with the semisimplicity of polarizable Hodge modules imply that $V_X^{\rm Hdg}$ is a direct factor of $\pi_*V_{X'}^{\rm Hdg}$. So we get a morphism
$$u:\iota_X(V_X^{\rm Hdg})\to\MM.$$
It is not clear whether $u$ is uniquely defined (since the decomposition is not unique). However, its underlying morphism of $\Q$-complexes coincides with the canonically defined adjunction morphism $j^\#$ so that  it   induces the desired morphism of mixed Hodge structures
$$ H^ij^\#: IH^i(X,\bV)\to H^i(U,\bV).$$

Let $\MM''$ be a mapping cone of $u:\iota_X(V_X^{\rm Hdg})\to\MM$ as defined in Section~\ref{sec4}. Remember \eqref{eqn:CrucialMHC} that $\MM$ comes from $j'_*V^{\rm Hdg}$, a mixed Hodge \emph{module}  of weight $\ge n+d$ (by Properties~\ref{Props}. F)) and hence $\Gr^W_k\MM=0$ for $k<d+n$. Then, by definition of the cone, one has 
\begin{equation} \label{eqn:lowgradeds}
\Gr^W_k\MM''=\Gr^W_k\MM=0\quad\hbox{for}\,\,\,k<d+n.
\end{equation}
Using Proposition~\ref{keyprop} (e.g.\ $\iota_X(V_X^{\rm Hdg})=\Gr^W_{d+n}H^0\MM$)  together with the long exact sequence~\eqref{eqn:longseq} we get moreover
$$\Gr^W_kH^i\MM''=0\quad\text{for}\,\,\,k\le i+d+n.$$
Since by \eqref{eqn:infty} we have $E_\infty^{i,k}= \Gr^W_kH^{i+k}\MM'$, the weight spectral sequence \eqref{eqn:specseq}  implies the surjectivity of
\[
\begin{diagram}
E_1^{-d-n-1, d+n+1+j} & \mapright{d_1} &  E_1^{-d-n,d+n+j+1}\\
\Vline{}{}{12pt}  && \Vline{}{}{12pt} \\
H^j\Gr^W_{d+n+1}\MM'' & \to &H^{j+1}\Gr^W_{d+n}\MM' 
\end{diagram}
\]
for all $j$, and this map splits by the semisimplicity of polarizable Hodge modules. So we get the surjectivity of
$$H^i(a_X)_*H^j\Gr^W_{d+n+1}\MM''\to H^i(a_X)_*H^{j+1}\Gr^W_{d+n}\MM''\quad\hbox{for any}\,\,i,j.$$
\begin{claim} This implies the surjectivity of
$$H^i(a_X)_*\Gr^W_{d+n+1}\MM''\to H^{i+1}(a_X)_*\Gr^W_{d+n}\MM''\quad\hbox{for any}\,\,i.$$
\end{claim}
\noindent\textit{Proof of the claim.} The truncation $\tau_{\le j}$ on $\Gr^W_kM''$ splits by the definition of mixed Hodge complexes so that  $\Gr^W_kM''\simeq \bigoplus_i  H^i(\Gr^W_kM'')[-i]$ where $M''$ is the underlying $D$-module of $\MM''$. Now  the truncation  induces a filtration $\tau'$ on $H^i(a_X)_*\Gr^W_kM''$ and  the preceding splitting for $\Gr^W_kM''$  coming from the truncation  induces a splitting  for $H^i(a_X)_*\Gr^W_{d+n+1}M''$ coming from $\tau'$. Its factors are isomorphic to $H^{i-j}(a_X)_*  H^j\Gr^W_{d+n+1}\MM''$ and this factor   maps surjectively to the factor  of $H^{i+1}(a_X)_*\Gr^W_{d+n}$ isomorphic to $H^{i-j}(a_X)_*H^{j+1}\Gr^W_{d+n}\MM''$. \qed
\\

\medskip
Again using $\Gr^W_k\MM''=0$ for $k<d+n$  \eqref{eqn:lowgradeds} it follows from the weight spectral sequence for $(a_X)_*\MM''$ that
\begin{equation}\label{eqn:final}
\Gr^W_kH^i(a_X)_*\MM''=0\quad\text{for}\,\,\,k\le d+n+i.
\end{equation}
{The   long exact sequence \eqref{eqn:longseq} for the direct image of the weakly distinguished triangle of the cone for $u$ under  $a_X:X\to $pt  reads
\[
\cdots  H^i(a_X)_*\MM'' \to  IH^i(X,\bV) \mapright{H^i(j^\#)}  H^i(U,\bV) \to H^{i+1}(a_X)_*\MM'' \, .
\]
From Corollary~\ref{stable}) this is a sequence of mixed Hodge structures and \eqref{eqn:final} shows the assertion about the lowest weights.}

{To complete the proof of Theorem~\ref{mainresult} in the analytic case, we  only have to show the independence of the compactification.  The map   \eqref{eqn:WhereItComesFrom} induced by $j$ is obtained from the natural map of $\Q$-complexes  $\iota: \ic_X(\bV[d])\to Rj_*\bV[d]$  after applying the global section functor. We may assume that we have a second  compactification $j':U\into Y$ related to $j:U\into X$ by a projective  morphism $f:X \to Y$.   The map induced by  $j'$ is similarly obtained from the natural map  $\iota': Rf_*\ic_X(\bV[d])\to Rf_*Rj_*\bV[d]=(Rj')_*\bV[d]$.
As before,  the decomposition Theorem~\ref{decomp} combined with the fact that  polarizable Hodge modules form a semi-simplicial category  implies  that the intersection complex $\ic_{Y}(\bV[d])$ is a direct factor of $Rf_*\ic_X(\bV[d])$  (but the latter might contain other direct factors). The restriction of $\iota'$ to  $\ic_{Y}(\bV[d])$  is exactly equal to $\iota$ while the other direct factors  are in the kernel of $\iota'$ since these must be supported on $Y \setminus U$. It follows that the image of $\iota$ does not depend on the compactification and hence neither does the image of  \eqref{eqn:WhereItComesFrom}.
 }

\end{document}